\documentclass{article}
\usepackage{amsthm,amsmath,amssymb}
\usepackage{algorithmic}
\usepackage{algorithm}
\usepackage{pst-tree,xspace}
\usepackage{array}
\usepackage{graphics}
\newtheorem{theorem}{Theorem}[section]
\newtheorem{proposition}{Proposition}[section]
\newtheorem{lemma}[theorem]{Lemma}
\theoremstyle{definition}
\newtheorem{definition}[theorem]{Definition}
\newtheorem{example}[theorem]{Example}

\theoremstyle{remark}

\numberwithin{equation}{section}

\newcommand{\binaryseptree}{binary separating tree\xspace}
\newcommand{\binaryseptrees}{binary separating trees\xspace}

\newcommand{\naryseptree}{contracted separating tree\xspace}
\newcommand{\naryseptrees}{contracted separating trees\xspace}

\newcommand{\intervaldecomptree}{common interval decomposition tree\xspace}
\newcommand{\Intervaldecomptree}{Common interval decomposition tree\xspace}
\newcommand{\intervaldecomptrees}{common interval decomposition trees\xspace}

\newcommand{\labeldecomptree}{labeled decomposition tree\xspace}
\newcommand{\Labeldecomptree}{Labeled decomposition tree\xspace}
\newcommand{\labeldecomptrees}{labeled decomposition trees\xspace}

\newcommand{\expanddecomptree}{expanded decomposition tree\xspace}
\newcommand{\expanddecomptrees}{expanded decomposition trees\xspace}

\begin{document}
\author{Mathilde Bouvel \thanks{LIAFA, Universit\'e Paris 7 and CNRS, Paris, France, 75251, mbouvel@liafa.jussieu.fr} \and Dominique Rossin\thanks{LIAFA, Universit\'e Paris 7 and CNRS, Paris, France, 75251, rossin@liafa.jussieu.fr, \texttt{http://www.liafa.jussieu.fr/$\sim$rossin}}}

\title{Longest Common Pattern between two Permutations} 

\maketitle 

\begin{abstract}
In this paper, we give a polynomial (${\mathcal O}(n^8)$)  algorithm for finding a longest common pattern between two permutations of size $n$ given that one is separable. We also give an algorithm for general permutations whose complexity depends on the length of the longest simple permutation involved in one of our permutations.

\end{abstract}

\section{Introduction and basic concepts}
The study of patterns in permutations has blossomed these last years: from a combinatorial point of view with the recent proof of the Stanley-Wilf conjecture by Marcus and Tard\"os, and from an algorithmic one with the development of algorithms for pattern and general pattern involvement. Although the general pattern involvement problem is $NP-$hard, some polynomial solutions exist for special kinds of patterns like the separable ones \cite{bbl,iba}. In this article we study the problem of finding a longest common pattern between two permutations $\sigma_1$ and $\sigma_2$ i.e. a permutation $\sigma$ which is involved in both permutations $\sigma_1$ and $\sigma_2$. This is a generalization of the {\em pattern involvement problem} since finding if the longest pattern between $\sigma_1$ and $\sigma_2$ is equal to $\sigma_1$ is equivalent to the pattern involvement problem.

First, we give a polynomial algorithm based on the work of \cite{bbl} for finding the longest common pattern if one permutation is separable. Then we generalize this algorithm for general permutations. The complexity of our algorithm is highly based on the length of the longest simple permutation \cite{bri} involved in our permutations.

\subsection{Permutations}

A permutation $\sigma$ of an interval $I$ of $\mathbb{N}$ is a bijective map from $I$ to itself. We denote by $\sigma_i$ the image of $i$ by $\sigma$. The permutation $\sigma$ could either be seen as a function or a word $\sigma_i \sigma_{i+1}\ldots\sigma_j$, where $I = \{h : i \leq h \leq j\}$. For example the permutation $\sigma=1 \ 4 \ 2 \ 5 \ 6 \ 3$ is the bijective function such that $\sigma(1)=1$, $\sigma(2)=4$, $\sigma(3)=2$, $\sigma(4)=5\ldots$. In the following, when we consider permutations without giving explicitly the interval $I$, we mean that $I = \{1, \ldots, n\}$ for some $n$.

\begin{definition}
  A permutation $\pi=\pi_1\ldots\pi_k$ is called a {\it pattern} of the permutation $\sigma=\sigma_{i+1}\ldots\sigma_{i+n}$ of $I = \{ h : i+1 \leq h \leq i+n\}$, with $k\leq n$,  if and only if there exist integers $i+1 \leq i_1 < i_2 <\ldots<i_k \leq i+n$ such that $\sigma_{i_{\ell}} < \sigma_{i_m}$ whenever $\pi_{\ell} < \pi_{m}$. We will also say that $\pi$ is involved in $\sigma$ or that $\sigma$ contains $\pi$. The subsequence $\sigma_{i_1} \sigma_{i_2}\ldots \sigma_{i_k}$ is called an \emph{occurrence} of $\pi$ in $\sigma$.
\end{definition}
A permutation $\sigma$ that does not contain $\pi$ as a pattern is said to {\em avoid} $\pi$. 
\begin{example}
For example $\sigma=1 \ 4 \ 2 \ 5 \ 6 \ 3$ contains the pattern $1 \ 3 \ 4 \ 2$ and $1 \ 5 \ 6 \ 3$, $1 \ 4 \ 6 \ 3$, $2\ 5\ 6\ 3$ and $1 \ 4 \ 5 \ 3$ are the occurrences of this pattern in $\sigma$. But $\sigma$ avoids the pattern $3 \ 2 \ 1$ as no subsequence of length $3$ of $\sigma$  is isomorphic to $3 \ 2 \ 1$ i.e. is decreasing.
\end{example}

A number of enumerative results has been proved on classes of pattern avoiding permutations for patterns of length $3,4$ and multiple patterns. More recently results about the algebricity of the generating function of general classes of permutations have been given \cite{bhv}.

Another field of study of these permutations is from the point of view of pattern involvement. The problem of deciding if a permutation $\pi$ is a pattern of a permutation $\sigma$ is $NP-$hard but this problem is proved to be polynomial if the pattern is separable \cite{bbl,iba}.

\begin{definition}
A permutation $\sigma$ of size $n$  is called {\em separable} if it avoids the patterns $3 \ 1 \ 4 \ 2$ and $2 \ 4 \ 1 \ 3$ or equivalently if it has a \binaryseptree.
\end{definition}

\begin{definition}
A {\em \binaryseptree} is a binary ordered tree with $n$ leaves such that each internal vertex is labeled by $+$ or $-$.
\par For each such tree, there is a unique way \cite{bbl} to decorate its leaves (considering them from left to right) by $\sigma_1,\sigma_2,\ldots,\sigma_n$ such that:
\begin{enumerate}
\item $\sigma_1\ldots\sigma_n$ is a permutation of $\{1 \ldots n\}$.
\item Each node (internal or leaf) is decorated by a permutation of an interval.
\item Each internal node labeled with $+$ (resp. $-$) is decorated by a permutation of $\{ i \ldots j\}$ such that its left child is decorated by a permutation of $\{i \ldots h-1\}$ (resp. $\{h\ldots j\}$) and its right child is decorated by a permutation of $\{h\ldots j\}$ (resp. $\{i \ldots h-1\}$), for some $h \in \{i+1, \ldots, j\}$.
\end{enumerate} \label{def-binaryseptree}
\end{definition}
It is easy to prove \cite{bbl} that to each separable permutation one can associate a \binaryseptree (see Figure \ref{fig:counterExampleSeparation}). Note that this tree is not uniquely defined as shown in figure \ref{fig:counterExampleSeparation}. However, one can associate a unique tree to each separable permutation by taking arbitrary ordered trees instead of binary ordered trees. These trees are the contraction of the binary ones by contracting every edge between two nodes with the same label $+$ or $-$. In these trees, the signs of the internal nodes are thus ranked: if the root of the tree has label $+$ (resp. $-$), then every node at odd depth has label $+$ (resp. $-$) and every node at even depth label $-$ (resp. $+$), so that all labels are determined from the label of the root.

\begin{definition}
The unique \emph{\naryseptree} associated to a separable permutation $\sigma$ is obtained from any \binaryseptree of $\sigma$ by contracting every edge between nodes of the same sign.\label{naryseptree}
\end{definition}

\begin{figure}
\begin{center}
\def\dedge{\ncline[linestyle=dashed,linecolor=red]}
\psset{levelsep=25pt,treesep=.5cm}
\pstree{\Tr{$+,423165897$}}
{
\pstree{\Tr[edge=\dedge]{$+,423165$}}
{
\pstree{\Tr{$-,4231$}}
{
\Tr{$4$}
\pstree{\Tr[edge=\dedge]{$-,231$}}
{
\pstree{\Tr{$+,23$}}{\Tr{$2$} \Tr{$3$}}
\Tr{$1$}
}
}
\pstree{\Tr{$-,65$}}{\Tr{$6$} \Tr{$5$}}
}
\pstree{\Tr{$-,897$}}
{
\pstree{\Tr{$+,89$}}{\Tr{$8$}\Tr{$9$}}
\Tr{$7$}
}
}
\pstree{\Tr{$+,423165897$}}
{
  \pstree{\Tr{$-,4231$}}{
    \pstree{\Tr[edge=\dedge]{$-,423$}} {
      \Tr{$4$}
      \pstree{\Tr{$+,23$}}{\Tr{$2$} \Tr{$3$}}
    }
    \Tr{$1$}
  }
  \pstree{\Tr[edge=\dedge]{$+,65897$}}{
    \pstree{\Tr{$-,65$}}{\Tr{$6$} \Tr{$5$}}
    \pstree{\Tr{$-,897$}}{
      \pstree{\Tr{$+,89$}}{\Tr{$8$}\Tr{$9$}}
      \Tr{$7$}
    }
  }
}
\rput[bl](6.2,-0.4){\rnode{Q}{label,decoration}}
\pstree{\Tr[name=R]{$+,423165897$}}
{
\pstree{\Tr{$-,4231$}}
{
\Tr{$4$}
\pstree{\Tr{$+,23$}}{\Tr{$2$} \Tr{$3$}}
\Tr{$1$}
}
\pstree{\Tr{$-,65$}}{\Tr{$6$} \Tr{$5$}}
\pstree{\Tr{$-,897$}}
{
\pstree{\Tr{$+,89$}}{\Tr{$8$}\Tr{$9$}}
\Tr{$7$}
}
}
\nccurve[angleB=90,angleA=0]{<-}{R}{Q}
\caption{Two decorated \binaryseptrees and the decorated \naryseptree of $\sigma=4 \ 2 \ 3 \ 1 \ 6 \ 5 \ 8 \ 9 \ 7$ }
\label{fig:counterExampleSeparation}
\end{center}
\end{figure}

\subsection{Modular decomposition of graphs ; Common interval decomposition of permutations \label{section-decomposition}}

\par
The \naryseptrees we introduced in Definition \ref{naryseptree} also appear in graph theory. Namely, those trees are a special case of \intervaldecomptrees (on which we however need to add a labeling). The \intervaldecomptrees are an equivalent for permutations of the modular decomposition trees for graphs \cite{fdm,bxhp}.
\par
Before we come to the pattern matching problem in permutations, we need to introduce the \intervaldecomptrees, the \labeldecomptrees, and finally the \expanddecomptrees that are the key structure we use in our algorithms.

\par
\begin{figure}[h!]
\begin{center}
\psset{unit=.2cm}
\begin{pspicture}(0,8)(32,-2)
\psset{linewidth=1.5pt}
\psline(0,0)(2,0)
\psline(3,0)(5,0)
\psline(6,0)(8,0)
\psline(9,0)(11,0)
\psline(12,0)(14,0)
\psline(15,0)(17,0)
\psline(18,0)(20,0)
\psline(21,0)(23,0)
\psline(24,0)(26,0)
\psline(27,0)(29,0)
\psline(30,0)(32,0)
\uput*{.2cm}[270]{0}(1,0){$5$}
\uput*{.2cm}[270]{0}(4,0){$1$}
\uput*{.2cm}[270]{0}(7,0){$10$}
\uput*{.2cm}[270]{0}(10,0){$9$}
\uput*{.2cm}[270]{0}(13,0){$6$}
\uput*{.2cm}[270]{0}(16,0){$7$}
\uput*{.2cm}[270]{0}(19,0){$8$}
\uput*{.2cm}[270]{0}(22,0){$11$}
\uput*{.2cm}[270]{0}(25,0){$2$}
\uput*{.2cm}[270]{0}(28,0){$4$}
\uput*{.2cm}[270]{0}(31,0){$3$}
\psline(27,1)(32,1)
\psline(24,2)(32,2)
\psline(12,3)(20,3)
\psline(6,5)(20,5)
\psline(6,6)(23,6)
\psline(0,7)(32,7)
\psset{linewidth=.5pt}
\psline(6,1)(11,1)
\psline(12,1)(17,1)
\psline(15,2)(20,2)
\psline(9,4)(20,4)
\end{pspicture}
%\vspace{2cm}\\
%\input{permdecomp.pstex_t}
\end{center}
\caption[Common interval decomposition of $\sigma = 5 \ 1 \ 10\ 9 \ 6 \ 7 \ 8 \ 11 \ 2 \ 4 \ 3$]{ Common interval decomposition of $\sigma = 5 \ 1 \ 10\ 9 \ 6 \ 7 \ 8 \ 11 \ 2 \ 4 \ 3$.\\ Common intervals are represented by horizontal lines, with strong intervals corresponding to bold lines. \label{decomposition}}
\end{figure}
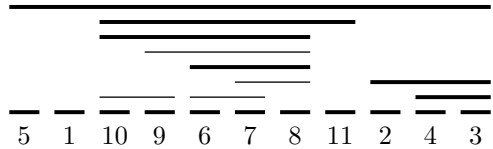

The \emph{common interval decomposition} of a permutation $\sigma$ of size $n$ (common meaning common with $Id_n$) is defined as follows. First consider all the \emph{common intervals} of $\sigma$ that is to say all the  subsequences $\sigma_j \sigma_{j+1} \ldots \sigma_k$ of consecutive entries of $\sigma$ such that $\{\sigma_j, \sigma_{j+1}, \ldots, \sigma_k\}$ is an interval of $\mathbb{N}$. Among the common intervals, the \emph{strong common intervals} (or strong intervals for short) are those that do not overlap any other common interval\footnote{The definition of overlapping common intervals follows the intuition: we say that two common intervals $I$ and $J$ are overlapping when $I \setminus J \neq \emptyset$, $J \setminus I \neq \emptyset$ and $I \cap J \neq \emptyset$}. Figure \ref{decomposition} illustrates the notion of common intervals and strong intervals. \\

\par 
The inclusion ordering yields a tree-like ordering on the set of strong intervals. 
This ordering is represented by a tree whose leaves are $\sigma_1$, $\sigma_2$, $\ldots$, $\sigma_n$ from left to right in this order, whose root is $\sigma$, and such that each internal node is the union of its children.

\par
Note that there are two different types of internal nodes in the tree of strong intervals. For some nodes, say $V$, with $k$ children $V_1,\ldots,V_k$ from left to right, there do not exist $(i,j) \neq (1,k)$ such that $1 \leq i < j \leq k$ and the union of $V_i,V_{i+1},\ldots,V_j$ is an interval. These nodes are called {\em prime} nodes and are of type $P$.

The other nodes $V$ are such that every union of consecutive children form an interval. Those node are called {\em linear} nodes and are of type $L$.

The tree along with the types $P$ and $L$ of the internal nodes (see Figure  \ref{arbredecomp}) is called the {\em \intervaldecomptree} of $\sigma$. In this tree, the order of the children of a node depends on $\sigma$ so that we have an ordered tree, unlike the modular decomposition trees for graphs.

Note that nodes of arity $2$ verify both linear and prime definitions. We choose to consider them of type $L$. This choice will be explained later. For further explanation on these trees, like the proof that a node is either linear or prime, refer to \cite{bxhp}. In \cite{bcmr,bxhp}, the authors show how to compute these trees in linear time.

\begin{figure}[h!]
\begin{center}
\pstree[nodesep=3pt,levelsep=*0.7,treesep=0.7]{\Tr[name=R]{{$P,5\, 1\, 10\, 9\, 6\, 7\, 8\, 11\, 2\, 4\,3$}}}{
	\begin{skiplevels}{3}
	     \Tr{5}
	     \end{skiplevels}
	\begin{skiplevels}{3}
	     \Tr{1}
	     \end{skiplevels}
	\pstree{\Tr{L,$10\, 9\, 6\, 7\, 8\, 11$}}{ 
		\pstree{\Tr{L,$10\, 9\, 6\, 7\, 8$}}{\begin{skiplevels}{1}\Tr{10} \Tr{9} \end{skiplevels} 
				\pstree{\Tr{L,$6\, 7\, 8$}}{\Tr{6} \Tr{7} \Tr{8} }}
		\begin{skiplevels}{2} \Tr{11} \end{skiplevels}}
	\begin{skiplevels}{1}
	 \pstree{\Tr{L,$2\, 4\, 3$}}
		{\begin{skiplevels}{1} \Tr{2} \end{skiplevels} \pstree{\Tr{L,$4\, 3$}}{\Tr{4} \Tr{3} }}
	\end{skiplevels}
}
\rput(-1,-0.5){\Rnode{P}{type,decoration}}
\nccurve[angleA=0,angleB=90]{<-}{R}{P}
\end{center}
\caption{\Intervaldecomptree for $\sigma = 5 \, 1 \, 10\, 9 \, 6 \, 7 \, 8 \, 11 \, 2 \, 4 \, 3$ \label{arbredecomp}}
\end{figure}

\par
For our algorithmic use of \intervaldecomptree, the types of the internal nodes are not sufficient and we need to label the internal nodes of the \intervaldecomptree of a permutation.
\par
It is easy to see that for each node $V$ of type $L$, the intervals of values corresponding to the children of $V$ are ordered either by increasing order or by decreasing order, when considering these children from left to right. The linear nodes are subsequently labeled $+$ or $-$ respectively. 
\par
The labeling of a  prime node with $d$ children consists of a permutation $\sigma$ of size $d$ that does not have any common intervals except $\{1\}, \ldots, \{d\}$ and $\sigma$. Such permutations are known as \emph{simple permutations} \cite{aa,aak,bri,bhv}. $\sigma$ represents the ordering of the children $V_1, \ldots , V_d$ of $V$ between them with respect to the values in the common intervals corresponding to the $V_i$'s. Namely, $\sigma_i < \sigma_j$ if and only if the common interval corresponding to $V_i$ contains values that are smaller than those contained in the common interval corresponding to $V_j$. For example, the simple permutation labeling the root of the tree on figure \ref{arbredecomp} is $3 \ 1 \ 4 \ 2$.

\par
This \intervaldecomptree along with labels $+$,$-$ and simple permutations can be computed easily in ${\mathcal O}(n^2 \ln n)$ which is sufficient for our purpose. Remind that the \intervaldecomptree can be computed in linear time. Thus it remains to add a label on each internal node $V$. This can be done by sorting the intervals corresponding to the children of $V$.

\par
When a \intervaldecomptree is labeled, the common intervals corresponding to the nodes can be deduced from this labeling, like in the case of \binaryseptrees (see definition \ref{def-binaryseptree}). So that the common intervals can be seen as a decoration on the nodes of the \intervaldecomptree with labels on internals nodes.
\begin{definition}
The \emph{\labeldecomptree} of a permutation $\sigma$ is the \intervaldecomptree of $\sigma$, where we add the labeling on internal nodes described above, and where we forget the decoration.
\end{definition}

\par
Figure \ref{labeldecomptree} gives the \labeldecomptree of $\sigma = 5 \ 1\  10\ 9 \ 6 \ 7 \ 8 \ 11 \ 2 \ 4 \ 3$.\\

\begin{figure}[h!]
\begin{center}
\pstree[nodesep=3pt,levelsep=*0.5,treesep=0.5]{\Tr[name=R1]{$3 \ 1 \ 4 \ 2$}} {
	\begin{skiplevels}{3}
	     \Tdot
	     \end{skiplevels}
	\begin{skiplevels}{3}
	     \Tdot
	     \end{skiplevels}
	\pstree{\Tr[name=R2]{$+$}}{ 
		\pstree{\Tr{$-$}}{\begin{skiplevels}{1}\Tdot  \Tdot \end{skiplevels} 
				\pstree{\Tr{$+$}}{\Tdot \Tdot \Tdot }}
		\begin{skiplevels}{2} \Tdot \end{skiplevels}}
	\begin{skiplevels}{1}
	 \pstree{\Tr{$+$}}
		{\begin{skiplevels}{1} \Tdot \end{skiplevels} \pstree{\Tr{$-$}}{\Tdot \Tdot}}
	\end{skiplevels}
}
\rput(0,0){\Rnode{P}{label}}
\nccurve[angleA=90,angleB=0]{->}{P}{R1}
\nccurve[angleA=180,angleB=0]{->}{P}{R2}
\pstree[nodesep=3pt,levelsep=*0.5,treesep=0.5]{\Tr{$3 \ 1 \ 4 \ 2$}} {
	\begin{skiplevels}{3}
	     \Tdot
	     \end{skiplevels}
	\begin{skiplevels}{3}
	     \Tdot
	     \end{skiplevels}
	\pstree{\Tr{$+$}}{ 
		\pstree{\Tr{$-$}}{\begin{skiplevels}{1}\pstree{\Tr{$-$}}{\Tdot  \Tdot} \end{skiplevels} 
				\pstree{\Tr{$+$}}{\pstree{\Tr{$+$}}{\Tdot \Tdot} \Tdot }}
		\begin{skiplevels}{2} \Tdot \end{skiplevels}}
	\begin{skiplevels}{1}
	 \pstree{\Tr{$+$}}
		{\begin{skiplevels}{1} \Tdot \end{skiplevels} \pstree{\Tr{$-$}}{\Tdot \Tdot}}
	\end{skiplevels}
}
\end{center}
\caption{\Labeldecomptree and \expanddecomptree for $\sigma = 5 \, 1 \, 10\, 9 \, 6 \, 7 \, 8 \, 11 \, 2 \, 4 \, 3$ \label{labeldecomptree}}
\end{figure}

\par
We have the following nice characterization of separable permutations in terms of \labeldecomptrees:
\begin{proposition}
The separable permutations are exactly those having a \intervaldecomptree with no prime nodes.\\
For any separable permutation, its \naryseptree and its \labeldecomptree are equal. \label{equivalence-septree-decomptree-1}
\end{proposition}

The proof of this proposition is straightforward from the definition of separable permutations. Note that it is important to give the type $L$ to binary internal nodes for stating that proposition.

\par
Proposition \ref{equivalence-septree-decomptree-1} states that \labeldecomptrees are a generalization to all permutations of the \naryseptree defined only for separable permutations. For \binaryseptrees, we have an analogous generalization, called {\em \expanddecomptrees}. In our algorithms, it is easier to work on binary nodes so that the natural representation used later is the \expanddecomptree.

\par
To transform a \labeldecomptree into an \expanddecomptree, take each linear node $V$ with children $V_1,\ldots,V_k$ and note that we can represent it by $(\ldots((V_1,V_2),V_3),\ldots),V_k)$ as shown in figure \ref{labeldecomptree}. Each  positive (resp. negative) internal node of arity $k$ gives $k-1$  positive (resp. negative) binary internal nodes.

A consequence of proposition \ref{equivalence-septree-decomptree-1} is:
\begin{proposition}
The \expanddecomptree of a separable permutation is one of its \binaryseptrees.
\label{equivalence-septree-decomptree}
\end{proposition}

\section{Longest common pattern between two permutations, including one separable}

\par
In this section, we describe a polynomial time algorithm for finding a longest common pattern between two permutations $\sigma$ and $\tau$ provided that $\sigma$ is separable. This algorithm uses the same technique as the one of Bose, Buss and Lubiw \cite{bbl} for finding an occurrence of a separable pattern in a general permutation. Namely, it computes a \binaryseptree $T_{\sigma}$ of $\sigma$ and uses it as a guide in the search of the longest common pattern with $\tau$. To avoid a complexity blow-up, the key point is to use dynamic programming, the initial problem being decomposed into sub-problems according to the structure of $T_{\sigma}$.

\par
First recall that it can be decided in linear time whether a permutation is separable or not. If it is, its \binaryseptree can also be computed in linear time, as described in \cite{bbl}. This result is actually  a special case of a more general one, proved in \cite{bcmr,bxhp}, and stating that the \intervaldecomptree of any permutation can be computed in linear time. 

\par
Instead of two permutations $\sigma$ and $\tau$, our algorithm takes as an input a \binaryseptree $T_{\sigma}$ for a separable permutation $\sigma$ of size $k$, and a permutation $\tau$ in the usual representation $\tau_1 \tau_2 \ldots \tau_n$. Notice that a \binaryseptree for $\sigma$ has $\mathcal{O}(k)$ nodes. \\

\par
More precisely, the algorithm fills in the array $$M=\{M(V,i,j,a,b) : V \textrm{ a node of } T_{\sigma}, 1 \leq i \leq j \leq n, 1 \leq a \leq b \leq n\}.$$ For any node $V$ in $T_{\sigma}$, let us denote by $\sigma(V)$ the subpermutation of $\sigma$ corresponding to the subtree of $T_{\sigma}$ rooted at $V$. With the notations of section 1, $\sigma(V)$ is the permutation decorating the node $V$. The cell $M(V,i,j,a,b)$ of the array $M$ contains a longest common pattern $\pi$ between $\sigma(V)$ on one hand, and the subpermutation $\tau_i \ldots \tau_j$ of $\tau$ on the other hand, with the additional restriction that the occurrence of $\pi$ in $\tau_i \ldots \tau_j$ uses only entries of $\tau$ whose values are between $a$ and $b$. The empty pattern, of size $0$, will be denoted $\epsilon$.
\begin{example}
If the node $V$ represents the pattern $2 \ 1$ (i.e. $\sigma(V) = (i+1)\ i$), and $\tau = 6 \ 4 \ 2 \ 5 \ 3 \ 1$ then we have for instance $M(V,2,4,3,5) = 1$, $M(V,2,5,3,4) = 2 \ 1$ and $M(V,4,5,1,2) = \epsilon$.
\end{example}
\par
The algorithm works as follows. To start the computation, we fill in the subarrays $M(V,\_,\_,\_,\_)$ for all the leaves $V$ of $T_{\sigma}$. Then, we compute $M(V,\_,\_,\_,\_)$ for any internal node $V$ using the subarrays $M(V_L,\_,\_,\_,\_)$ and $M(V_R,\_,\_,\_,\_)$ corresponding to the left child ($V_L$) and the right child ($V_R$) of $V$. In order to combine the patterns found in $M(V_L,\_,\_,\_,\_)$ and $M(V_R,\_,\_,\_,\_)$ with the intention of filling $M(V,\_,\_,\_,\_)$, we need a definition of pattern concatenation.
\begin{definition}
Consider two patterns $\pi$ and $\pi'$ of respective lengths $k$ and $k'$. The positive and negative concatenations of $\pi$ and $\pi'$ are defined respectively by:
$\pi \oplus \pi' = \pi_1 \cdots \pi_k (\pi'_1 +k) \cdots (\pi'_{k'} +k) \textrm{ and } \pi \ominus \pi' = (\pi_1+k') \cdots (\pi_k+k') \pi'_1 \cdots \pi'_{k'} $ \label{concatenation}
\end{definition}
\begin{example}~\\
$4 \ 3 \ 5 \ 2 \ 1 \oplus 3 \ 1 \ 4 \ 2 = 4 \ 3 \ 5 \ 2 \ 1 \ \vdots \ 8 \ 6 \ 9 \ 7$ \\
$ 4 \ 3 \ 5 \ 2 \ 1 \ominus 3 \ 1 \ 4 \ 2= 8 \ 7 \ 9 \ 6 \ 5 \ \vdots \ 3 \ 1 \ 4 \ 2$
\end{example}

It is clear from definition \ref{concatenation} that the positive (resp. negative) concatenation of two patterns produces again a pattern.\\

\par
The detailed dynamic programming algorithm is given in Algorithm \ref{alg:separable}.

\begin{algorithm}[h!t]
\caption{Longest common pattern between two permutations, first one separable}
\label{alg:separable}
\begin{algorithmic}[1]
\STATE {\textsc{Input}: A \binaryseptree $T_\sigma$ of a separable permutation $\sigma$ of size $k$ and a permutation $\tau$ of size $n$ \\~\\}
\STATE {\textsc{Create an array $M$:} }
\FOR{any node $V$ of $T_\sigma$ and any integers $i$, $j$, $a$ and $b$ between $1$ and $n$}
\STATE {$M(V,i,j,a,b) \leftarrow \epsilon$ }
\ENDFOR \\~\\
\STATE {\textsc{Fill in $M$ for the leaves of $T_\sigma$:} } 
\FOR{any leaf $V$ of $T_\sigma$}
\FOR{any integers $i$, $j$, $a$ and $b$ between $1$ and $n$, $i\leq j$, $a \leq b$}
\IF {there exists some $h \in \{i, i+1, \ldots, j \}$ such that $a \leq \tau_{h} \leq b$}
\STATE {$M(V,i,j,a,b) \leftarrow 1$ }
\ENDIF
\ENDFOR
\ENDFOR \\~\\
\STATE {\textsc{Fill in the rest of $M$:} }
\FOR{any internal node $V$ of $T_\sigma$, considering the nodes in the postfix ordering} 
\IF {$V$ is a positive node}
\FOR{any integers $i$, $j$, $a$ and $b$ between $1$ and $n$, $i\leq j$, $a \leq b$}
\STATE {$M(V,i,j,a,b) \leftarrow Longest(\{  M(V_L,i,h-1,a,c-1) \oplus M(V_R,h,j,c,b)  :$ \\ \qquad \qquad \qquad \qquad \qquad \qquad $ i \leq h \leq j+1, a \leq c \leq b+1\}) $
}
\ENDFOR
\ELSE \STATE { /* $V$ is a negative node /*}
\FOR{any integers $i$, $j$, $a$ and $b$ between $1$ and $n$, $i\leq j$, $a \leq b$}
\STATE {$M(V,i,j,a,b) \leftarrow Longest(\{  M(V_L,i,h-1,c,b) \ominus M(V_R,h,j,a,c-1)  :$ \\ \qquad \qquad \qquad \qquad \qquad \qquad $ i \leq h \leq j+1, a \leq c \leq b+1\}) $}
\ENDFOR
\ENDIF
\ENDFOR \\~\\
\STATE {\textsc{Output}: $M(\textrm{root of }T_\sigma,1,n,1,n)$}
\end{algorithmic}
\end{algorithm}

In the previous, for any set $S$ of patterns, $Longest(S)$ returns a longest pattern in the set $S$. If we want $Longest$ to be a uniquely defined application, we can choose $Longest(S)$ as the smallest pattern with respect to the lexicographical ordering among the patterns in $S$ having maximal length.

\begin{proposition}
Algorithm \ref{alg:separable} is correct: it outputs a longest common pattern between the two permutations $\sigma$ and $\tau$ given in input. \label{correctness}
\end{proposition}

\begin{proof}
The proof is by induction. \\ We show that the algorithm stores in $M(V,i,j,a,b)$ a longest common pattern between $\sigma(V)$ and $\tau_i \ldots \tau_j$ whose occurrence in $\tau_i \ldots \tau_j$ uses only values between $a$ and $b$.

\par
If $V$ is a leaf, the above statement is clearly true.

\par
If $V$ is an internal node, with two children $V_L$ (its left child) and $V_R$ (its right child), then let $i$, $j$, $a$ and $b$ be integers such that $1 \leq i \leq j \leq n$ and $1 \leq a \leq b \leq n$. We assume in the rest of the proof that $V$ is a positive node, the case of a negative node being very similar. To begin with, it is easy to see that $M(V,i,j,a,b)$ contains a common pattern between $\sigma(V)$ and $\tau_i \ldots \tau_j$ using only values between $a$ and $b$ in $\tau_i \ldots \tau_j$. Indeed, by induction hypothesis, we infer that every pattern in the set $S = \{  M(V_L,i,h-1,a,c-1) \oplus M(V_R,h,j,c,b)  : i \leq h \leq j+1, a \leq c \leq b+1\}$, and \emph{a fortiori} $Longest(S)$, is a common pattern between $\sigma(V)$ and $\tau_i \ldots \tau_j$ using only values between $a$ and $b$ in $\tau_i \ldots \tau_j$.
\begin{example} Consider the permutations $\sigma=1 \ 4 \ 2 \ 3 \ 6 \ 5 \ 7 \ 8$ and \\
\begin{tabular}{p{7.5cm}c}
$\tau =4 \ 1 \ 3 \ 2  \ 5 \ 6 \ 8 \ 9 \ 7$. A separating tree $T_{\sigma}$ of $\sigma$ is represented on the right. For example, we can choose $V$ to be the right child of the root of $T_{\sigma}$. Then we have $\sigma(V)= 4 \ 2 \ 3 \ 6 \ 5 \ 7 \ 8$, $\sigma(V_L) = 4 \ 2 \ 3$ and $\sigma(V_R) = 6 \ 5 \ 7 \ 8$. Now choose $i = 2, j = 7, a =2$ and $b=8$. We want to show that for any $h \in \{i, \ldots j+1\}$ and $c \in \{a, \ldots b+1\}$, $M(V_L,i,h-1,a,c-1) \oplus M(V_R,h,j,c,b)$ is a common pattern between $\sigma(V)$ and $\tau_i \ldots \tau_j$ using only values between $a$ and $b$ in $\tau_i \ldots \tau_j$. Take for example $h = 5$ and $c =4$.
& \pstree[nodesep=3pt,levelsep=*0.7,treesep=0.3]{\Tr{+}}
{\begin{skiplevels}{3}\Tr{1}\end{skiplevels} \pstree{\Tr{+}}{ 
\pstree{\Tr{-}}{ \skiplevel{\Tr{4}} \pstree{\Tr{+}}{ \Tr{2} \Tr{3}}} 
\pstree{\Tr{+}}{\pstree{\Tr{-}}{\Tr{6} \Tr{5}} \pstree{\Tr{+}}{\Tr{7} \Tr{8}}} }
}
\end{tabular}\\
By induction hypothesis, $M(V_L,i,h-1,a,c-1)$ contains a longest common pattern between $\sigma(V_L)$ and $\tau_i \ldots \tau_{h-1}$ using only values between $a$ and $c-1$ in $\tau_i \ldots \tau_{h-1}$. Here, $M(V_L,i,h-1,a,c-1) = 2 \ 1$, an occurrence of $2 \ 1$ in $\tau_i \ldots \tau_{h-1} = 1 \ 3 \ 2$ using values between $2$ and $3$ being $1 \ \textbf{3} \ \textbf{2}$, and an occurrence of $2 \ 1$ in $\sigma(V_L) = 4 \ 2 \ 3$ being $ \textbf{4} \ \textbf{2} \ 3$. Similarly, we have $M(V_R,h,j,c,b) =M(V_R, 5, 7, 4, 8) = 1 \ 2 \ 3$, as shown by the occurrences $\textbf{5} \ \textbf{6} \ \textbf{8}$ in $\tau_h \ldots \tau_j=5 \ 6 \ 8$ and $\textbf{6} \ 5 \ \textbf{7} \ \textbf{8}$ in $\sigma(V_R) = 6 \ 5 \ 7 \ 8$. The occurrence of $M(V_L,i,h-1,a,c-1) \oplus M(V_R,h,j,c,b) = 2 \ 1 \ 3 \ 4 \ 5$ in $\tau_i \ldots \tau_j$ using values between $a$ and $b$ is thus obtained by considering simultaneously the two occurrences in $\tau$ enlightened before. Namely, an occurrence of $2 \ 1 \ 3 \ 4 \ 5$ in $\tau_i \ldots \tau_j = 1 \ 3 \ 2 \ 5 \ 6 \ 8$ using values between $2$ and $8$ is $1 \ \textbf{3} \ \textbf{2} \ \textbf{5} \ \textbf{6} \ \textbf{8}$. Notice also that the occurrence $\textbf{4} \ \textbf{2} \ 3 \ \textbf{6} \ 5 \ \textbf{7} \ \textbf{8}$ of $2 \ 1 \ 3 \ 4 \ 5$ in $\sigma(V)= 4 \ 2 \ 3 \ 6 \ 5 \ 7 \ 8$ is again obtained considering simultaneously the occurrences of $2 \ 1$ and $1 \ 2 \ 3$ in $\sigma(V_L)$ and $\sigma(V_R)$ respectively.
\end{example}

\par
It remains to show that:
\begin{lemma}
\label{lem:longest}
$Longest(S)$ is of maximal length among all the common patterns between $\sigma(V)$ and $\tau_i \ldots \tau_j$ using only values between $a$ and $b$ in $\tau_i \ldots \tau_j$.
\end{lemma}
\begin{proof}
\par
We denote by $\pi$ a longest common pattern between $\sigma(V)$ and $\tau_i \ldots \tau_j$, using only values between $a$ and $b$ in $\tau_i \ldots \tau_j$.\\
As shown on Figure \ref{fig:proofLemma}, there exist integers $h \in \{i, \ldots , j+1\}$ and $c \in \{a, \ldots, b+1\}$ such that $\pi$ is decomposed into $\pi = \pi_1 \oplus \pi_2$, with $\pi_1$ a common pattern between $\sigma(V_L)$ and $\tau_i \ldots \tau_{h-1}$, using only values between $a$ and $c-1$ in $\tau_i \ldots \tau_{h-1}$, and $\pi_2$ a common pattern between $\sigma(V_R)$ and $\tau_h \ldots \tau_j$, using only values between $c$ and $b$ in $\tau_h \ldots \tau_j$. \\

\begin{figure}[ht]
\begin{center}
\input{dessin.pstex_t}
\caption{Proof of Lemma \ref{lem:longest}}
\label{fig:proofLemma}
\end{center}
\end{figure}

Notice that in this decomposition $\pi_1$ or $\pi_2$ might be the empty pattern. \\
It can be easily seen that if $\pi_1$ (respectively $\pi_2$) were not a pattern of maximal length for the given intervals of indices and values, then $\pi$ would not be of maximal length either, contradicting the definition of $\pi$. Consequently, $\pi_1$ and $\pi_2$ are longest common patterns between $\sigma(V_L)$ and $\sigma(V_R)$ respectively, and $\tau$ in the given intervals of indices and values. By induction hypothesis, $|M(V_L,i,h-1,a,c-1)| = | \pi_1 |$ and $|M(V_R,h,j,c,b)| = | \pi_2 |$. The pattern stored in $M(V,i,j,a,b)$ by the algorithm is of length at least $|M(V_L,i,h-1,a,c-1) \oplus M(V_R,h,j,c,b)| = | \pi_1 \oplus \pi_2 | = | \pi |$. $\pi$ being of maximal length by assumption, we conclude that $M(V,i,j,a,b)$ is also of maximal length.
\end{proof}
Finally, $M(V,i,j,a,b)$ contains a longest common pattern between $\sigma(V)$ and $\tau_i \ldots \tau_j$ whose occurrence in $\tau_i \ldots \tau_j$ uses only values between $a$ and $b$.
\par
When $V$ is a negative node, we decompose $\pi$ into $\pi_1 \ominus \pi_2$, with $\pi_1$ a common pattern between $\sigma(V_L)$ and $\tau_i \ldots \tau_{h-1}$ , using only values between $c$ and $b$ in $\tau_i \ldots \tau_{h-1}$, and $\pi_2$ a common pattern between $\sigma(V_R)$ and $\tau_h \ldots \tau_j$, using only values between $a$ and $c-1$ in $\tau_h \ldots \tau_j$, and the proof follows the same steps.
\end{proof}

\begin{proposition}
Algorithm \ref{alg:separable} has a time complexity in $\mathcal{O}(min(k,n) k n^6)$. \label{polynomial}
\end{proposition}

\begin{proof}
Algorithm \ref{alg:separable} handles an array $M$ of size $\mathcal{O}(k n^4)$, where each cell contains a pattern of length at most $min(n,k)$, so that the total space complexity is $\mathcal{O}(min(n,k) k n^4)$. For filling in the subarrays $M(V,\verb=_=,\verb=_=,\verb=_=,\verb=_=)$ for all the leaves $V$ of $T_{\sigma}$ (line 6 to 13 of algorithm \ref{alg:separable}), the total time complexity is $\mathcal{O}(k n^5)$. And for any internal node $V$, filling in one entry of the subarray $M(V,\verb=_=,\verb=_=,\verb=_=,\verb=_=)$ costs $\mathcal{O}(min(k,n) n^2 )$, since at line 18 (or 23) of algorithm \ref{alg:separable}, we search for an element of maximal length among $\mathcal{O}( n^2 )$ elements, each of size $\mathcal{O}(min(k,n))$. Consequently, filling in completely this subarray $M(V,\verb=_=,\verb=_=,\verb=_=,\verb=_=)$ requires a time complexity $\mathcal{O}(min(k,n) n^6 )$. Since there are $\mathcal{O}(k)$ internal nodes in $T_{\sigma}$, we have the announced result. 
\end{proof}

This complexity can be improved to $\mathcal{O}(k n^6)$, storing an integer, a label ($\oplus$ or $\ominus$), and two pointers in $M(V,i,j,a,b)$ (when $V$ is an internal node) instead of a pattern. Namely, if algorithm \ref{alg:separable} fills in $M(V,i,j,a,b)$ with the pattern $\rho = M(V_L,i,h-1,a,c-1) \oplus M(V_R,h,j,c,b)$, it is sufficient to store in $M(V,i,j,a,b)$ the length of $\rho$, the label $\oplus$, and two pointers pointing to the entries $M(V_L,i,h-1,a,c-1)$ and $M(V_R,h,j,c,b)$ of the array $M$. At the end of the algorithm, this system of pointers gives a \binaryseptree of a longest common pattern $\pi$ between $\sigma$ and $\tau$. From this tree, $\pi$ can be computed in linear time \cite{bbl}.\\

\par
A consequence of properties \ref{correctness} and \ref{polynomial} is: 

\begin{theorem}
 The problem of finding a longest common pattern between two permutations, the first one being separable, is in $P$. \label{thm-separable}
\end{theorem}

\section{Longest common pattern between two permutations}

\par
The result of theorem \ref{thm-separable} can be easily extended to classes of permutations that are less restricted than separable permutations. Using the \intervaldecomptree introduced in section \ref{section-decomposition}, we will see that a longest common pattern between two permutations $\sigma$ and $\tau$ can be computed in polynomial time as soon as the arity of any prime node in the \intervaldecomptree of $\sigma$ is bounded by a constant $d$ chosen independently.\\

\par
Proposition \ref{equivalence-septree-decomptree} states that \expanddecomptrees are a generalization to all permutations of \binaryseptrees, defined only for the separable ones. From this remark, it becomes natural to try and use \expanddecomptrees in an algorithm for finding a longest common pattern between two general permutations. In the following, we describe such an algorithm and analyze its complexity: it is not a polynomial time algorithm, but the complexity analysis reveals classes of permutations for which the algorithm runs in polynomial time.\\

\par
First, we notice that the \expanddecomptree of a permutation can be computed in polynomial time, and even in linear time if we have an independent bound on the arity of the prime nodes. \cite{bcmr,bxhp} provide a $\mathcal{O}(n)$-time algorithm for computing the \intervaldecomptree $T_{\sigma}$ of any permutation $\sigma$ of size $n$. This tree can easily be labeled in time $\mathcal{O}(d \log d \cdot n)$ where $d$ is a bound on the arity of prime nodes: for every internal node $V$ of $T_{\sigma}$ -- that are in number $\mathcal{O}(n)$ -- it is enough to sort its children to find the label of $V$. Finally, the vertical expansion necessary to obtain the \expanddecomptree of $\sigma$ requires again a linear time. \\

\par
Algorithm \ref{alg:decomp} takes as an input an \expanddecomptree $T_{\sigma}$ of a permutation $\sigma$ and a permutation $\tau$. It outputs a longest common pattern between $\sigma$ and $\tau$. It works just like algorithm \ref{alg:separable}, except for the case of prime nodes in $T_{\sigma}$. The procedure in this additional case is described in algorithm \ref{alg:decomp}.\\

\begin{algorithm}[!ht]
\caption{Longest common pattern between two permutations}
\label{alg:decomp}
\begin{algorithmic}[1]
\STATE {\textsc{Input} : An \expanddecomptree $T_\sigma$ of a permutation $\sigma$ of size $k$ and a permutation $\tau$ of size $n$ \\~\\}
\STATE {\textsc{Create an array $M$}: proceed as in algorithm \ref{alg:separable} \\~\\}
\STATE {{\sc Fill in $M$ for the leaves of $T_\sigma$}: proceed as in algorithm \ref{alg:separable}\\~\\}
\STATE {{\sc Fill in the rest of $M$}: }
\FOR{any internal node $V$ of $T_\sigma$, considering the nodes in the postfix ordering} 
\IF {$V$ is a positive or a negative node}
\STATE {proceed as in algorithm \ref{alg:separable}}
\ELSE \STATE{ /* $V$ is a prime node /*}
\STATE {Let $\rho$ be the simple permutation labeling $V$}
\STATE{Let $d$ be the arity of $V$, and $V_1, \ldots , V_d$ the children of $V$, from left to right}
\FOR{any integers $i$, $j$, $a$ and $b$ between $1$ and $n$, $i\leq j$, $a \leq b$}
\STATE {$M(V,i,j,a,b) \leftarrow Longest(S) \textrm{ where }$ \\ $ S =\big\{  \odot_{\rho}(\big( M(V_k, h_{k-1}, h_k-1, c_{\rho_k -1}, c_{\rho_k} -1) \big)_{1 \leq k \leq d} )  : $ \\ $\qquad  i = h_0 \leq h_1 \leq \ldots \leq h_d = j+1, a = c_0 \leq c_1 \leq \ldots \leq h_d = b+1 \big\} $}
\ENDFOR
\ENDIF
\ENDFOR \\~\\
\STATE {\textsc{Output} : $M(\textrm{root of }T_\sigma,1,n,1,n)$}
\end{algorithmic}
\end{algorithm} 

\par
Algorithm \ref{alg:decomp} uses a more general kind of pattern concatenation than just $\oplus$ and $\ominus$. The $\rho$-concatenation, or pattern concatenation according to $\rho$, is defined as follows:
\begin{definition}
Given $\rho = \rho_1 \ldots \rho_n$ a permutation of size $n$, and $n$ patterns $\pi^1, \ldots, \pi^n$ of respective size $k_1, \ldots , k_n$, the $\rho$-concatenation of the $(\pi^i)_{1 \leq i \leq n}$ is 
$$\odot_{\rho}(\pi^1, \ldots, \pi^n) = \textrm{shift}(\pi^1,\rho_1) \ldots \textrm{shift}(\pi^n,\rho_n) \textrm{ where } $$
$$
\textrm{shift}(\pi^i,\rho_i) = \textrm{shift}(\pi^i,\rho_i)(1) \ldots \textrm{shift}(\pi^i,\rho_i)(k_i) \textrm{ and } 
$$
$$
\textrm{shift}(\pi^i,\rho_i) (x) = (\pi^i(x) + k_{\rho^{-1}_1} + \ldots + k_{\rho^{-1}_{i-1}}) \textrm{ for all } x \textrm{ between } 1 \textrm{ and } k_i
$$
\end{definition}

\begin{example}~\\
$\odot_{25314}(21,312,4321,12,231) = 4 \ 3 \ \vdots \ 14 \ 12 \ 13 \ \vdots \ 8 \ 7 \ 6 \ 5 \ \vdots \ 1 \ 2 \ \vdots \ 10 \ 11 \ 9$.
\end{example}

We can also notice that $\odot_{12} = \oplus$ and $\odot_{21} = \ominus$.\\

\par
The idea behind algorithm \ref{alg:decomp} is quite simple. When filling in $M(V,i,j,a,b)$ for a prime node $V$ labeled by $\rho$ and having $d$ children $V_1 \ldots V_d$, we ``slice'' the intervals $\{i, \ldots , j\}$ and $\{a,\ldots,b\}$ into $I_1 \ldots I_d$ and $A_1 \ldots A_d$ respectively, such that $I_p \prec I_k$ and $A_p \prec A_k$\footnote{By $A\prec B$, we mean that $\forall a \in A, \forall b \in B, a < b$.} whenever $p < k$. Then we $\rho$-concatenate longest common patterns between the $\sigma(V_k)$ and $\tau$ in the intervals $I_k$ of indices and $A_{\rho_k}$ of values. With the notation of algorithm \ref{alg:decomp}, $I_k = \{h_{k-1}, \ldots, h_k-1\}$ and $A_k = \{c_{k-1}, \ldots, c_k-1\}$.

\begin{proposition}
Algorithm \ref{alg:decomp} is correct: it outputs a longest common pattern between the two permutations $\sigma$ and $\tau$ given in input. \label{correctness2}
\end{proposition}

\begin{proof}
Similar to the proof of proposition \ref{correctness}.
\par
With the notations of the proof of proposition \ref{correctness}, in the case of a prime node $V$ labeled by $\rho$, with children $V_1 \ldots V_d$, there exist integers $i = h_0 \leq h_1 \leq \ldots \leq h_d = j+1$ and $a = c_0 \leq c_1 \leq \ldots \leq c_d = b+1$, such that we can decompose $\pi$ into $\pi = \odot_{\rho}(\pi^1, \ldots, \pi^d)$, with $\pi^k$ a common pattern between $\sigma(V_k)$ and $\tau_{h_{k-1}} \ldots \tau_{h_k -1}$ , using only values between $c_{\rho_k -1}$ and $c_{\rho_k}-1$ in $\tau_{h_{k-1}} \ldots \tau_{h_k -1}$. Using this decomposition of $\pi$, we can use the induction hypothesis on the nodes $(V_k)_{1\leq k \leq d}$ and finish the proof as before.
\end{proof}

\par
In this proof, the trick relies on the fact that a common pattern between $\sigma(V)$ and $\tau_i \ldots \tau_j$ is always a concatenation of common patterns between the children of $V$ and ``slices'' of $\tau_i \ldots \tau_j$. This stability when going from parents to children in the \expanddecomptree also appears in a paper of Albert and Atkinson \cite{aa}, for example in their lemma 15. \\

\par
The main difference between algorithms \ref{alg:separable} and  \ref{alg:decomp} lies in the complexity analysis. Those two algorithms deal with dynamic programming arrays of the same size, but the cost for computing one entry can be very superior in algorithm \ref{alg:decomp} than in algorithm \ref{alg:separable}. Indeed, for any internal node $V$, in order to fill in one entry of $M(V,\verb=_=,\verb=_=,\verb=_=,\verb=_=)$, algorithm \ref{alg:separable} computes a longest pattern in a set containing $\mathcal{O}(n^2)$ elements, whereas in algorithm \ref{alg:decomp}, the set from which we have to extract a longest pattern contains $\mathcal{O}(n^{2d-2})$ elements, if $d$ is the arity of $V$ (see line 13 of algorithm \ref{alg:decomp}).

\par
With no hypothesis on a permutation $\sigma$ of size $k$, the only bound we can give on the maximal arity $d$ of a prime node in the \expanddecomptree of $\sigma$ is $d \leq k$. This bound is optimal since the equality $d=k$ is achieved when $\sigma$ is a simple permutation. The total time complexity of algorithm  \ref{alg:decomp} is consequently $\mathcal{O}(min(n,k)kn^{2k+2})$, and it is not polynomial. However, if we consider classes of permutations such that the arity of any prime node in their \expanddecomptree is bounded by a constant $d$, algorithm \ref{alg:decomp} has a time complexity $\mathcal{O}(min(n,k)kn^{2d+2})$. In this particular case, algorithm \ref{alg:decomp} runs in polynomial time.

\par
This can be summarized in the following theorem:
\begin{theorem}
Let $d$ be a integer. Consider the class $R$ of permutations having an \expanddecomptree with all prime nodes of arity smaller than $d$. Then the problem of finding a longest common pattern between a permutation in $R$ and another unrestricted permutation is in $P$. \label{thmfinal}
\end{theorem}

\section{Conclusion and open problems}

We generalize the algorithm given in \cite{bbl} for the problem {\em longest common pattern}. Yet our algorithm seems far from optimal. For example, for separable permutations, our work is based on \cite{bbl} but Ibarra \cite{iba} give a faster $({\mathcal O}(n^5))$ algorithm for the pattern involvement problem. Could this algorithm be adapted to the {\em longest common pattern} problem? Yet a lower bound is given by the edit distance problem \cite{zhangshasha} as the edit distance problem between two trees is a special case of the longest common pattern problem as shown in \cite{micheliRossin}.

%use either amsalpha bibliography style or acm

%use either amsalpha bibliography style or acm

\bibliographystyle{amsalpha}
\bibliography{Biblio}
% \bibitem{Bour}
% N. Bourbaki,
%   Groupes et Alg\`ebres de Lie, Chapitres 4, 5 et 6, Hermann,
%   Paris (1968).

\end{document}